\documentclass[12pt]{article}
\usepackage{amsmath,amsthm,amssymb}

\newtheorem{prop}{Proposition}

\newtheorem{rem}{Remark}

\begin{document}

\newcommand{\Rg}       {\mathbb{R}}
\newcommand{\Pg}       {\mathbb{P}}
\newcommand{\Eg}       {\mathbb{E}}
\newcommand{\ungras}{\mathbb{1}}

\newcommand{\sRR}{{\hbox{$\scriptstyle{I}$\kern-.25em\hbox{$\scriptstyle R$}}}}

\newcommand{\equiva}    {\displaystyle\mathop{\simeq}}
\newcommand{\eqdef}     {\stackrel{\triangle}{=}}
\newcommand{\limps}     {\mathop{\hbox{\rm lim--p.s.}}}
\newcommand{\Limsup}    {\mathop{\overline{\rm lim}}}
\newcommand{\Liminf}    {\mathop{\underline{\rm lim}}}
\newcommand{\vers}      {\mathop{\;{\rightarrow}\;}}
\newcommand{\abs}[1]{\left| #1 \right|}

\newcommand{\AAA}{{\cal A}}
\newcommand{\BB}{{\cal B}}
\newcommand{\CC}{{\cal C}}
\newcommand{\DD}{{\cal D}}
\newcommand{\EE}{{\cal E}}
\newcommand{\FF}{{\cal F}}
\newcommand{\GG}{{\cal G}}
\newcommand{\HH}{{\cal H}}
\newcommand{\II}{{\cal I}}
\newcommand{\JJ}{{\cal J}}
\newcommand{\KK}{{\cal K}}
\newcommand{\LL}{{\cal L}}
\newcommand{\NN}{{\cal N}}
\newcommand{\MM}{{\cal M}}
\newcommand{\OO}{{\cal O}}
\newcommand{\PP}{{\cal P}}
\newcommand{\QQ}{{\cal Q}}
\newcommand{\RR}{{\cal R}}
\newcommand{\SSS}{{\cal S}}
\newcommand{\TT}{{\cal T}}
\newcommand{\UU}{{\cal U}}
\newcommand{\VV}{{\cal V}}
\newcommand{\WW}{{\cal W}}
\newcommand{\XX}{{\cal X}}
\newcommand{\YY}{{\cal Y}}
\newcommand{\ZZ}{{\cal Z}}
\newcommand{\carre}{\hfill\rule{0.25cm}{0.25cm}}
\newcommand{\tr}   {\mbox{ \rm{tr} }}
\newcommand{\noi}  {\noindent}
\newcommand{\non}  {\nonumber}
\renewcommand{\EE}   {{\bf E}}

\renewcommand{\theenumi}{\roman{enumi}}
\renewcommand{\labelenumi}{{\rm({\it\theenumi}\/)}}
\renewcommand{\theenumii}{\alph{enumii}}
\renewcommand{\labelenumii}{{\it\theenumii.}}

\newcommand{\bbN}{\mathbb{N}}
\newcommand{\bbZ}{\mathbb{Z}}
\newcommand{\bbQ}{\mathbb{Q}}
\newcommand{\bbR}{\mathbb{R}}
\newcommand{\bbC}{\mathbb{C}}

\newcommand{\bbP}{\mathbb{P}}

\newcommand{\mcF}{\mathcal{F}}

\newcommand{\eps}{\varepsilon}

\newcommand{\const}{\mathrm{const}}

\newcommand{\wtW}{\widetilde{W}}
\newcommand{\wh}{\widehat{h}}
\newcommand{\bg}{\bar{\gamma}}

\newcommand{\argmin}{\mathop{\mathrm{argmin}}}
\renewcommand{\tanh}{\mathop{\mathrm{th}}}
\renewcommand{\sinh}{\mathop{\mathrm{sh}}}
\renewcommand{\cosh}{\mathop{\mathrm{ch}}}

\newcommand{\fexp}[1]{\exp\left\{ #1 \right\}}




\pagenumbering{arabic}

\title{About Gaussian filtering problems with general exponential quadratic  criteria}
\author{
M.L.~Kleptsyna \footnote{Author for correspondence.
Fax: +33 2 43 83 35 79.}\\
\vspace{-2mm} \small\sl
Laboratoire de Statistique et Processus / Universit\'e du Maine\\
\vspace{-2mm}
\small\sl Av. Olivier Messiaen, 72085 Le Mans, Cedex 9, France \\
\small\sl e-mail: Marina.Kleptsyna@univ-lemans.fr\\
\mbox{\hspace{100mm}}\vspace{-5mm} \and \protect{\newline}
\vspace{-1mm}
A.~Le Breton \\
\vspace{-2mm} \small\sl
Laboratoire Jean Kuntzmann / Universit\'e J. Fourier\\
\vspace{-2mm}
\small\sl BP 53, 38041 Grenoble Cedex 9, France  \\
\small\sl e-mail: Alain.Le-Breton@imag.fr\\
\mbox{\hspace{100mm}}\vspace{-5mm} \and \protect{\newline}
\vspace{-1mm} M.~Viot\\
\vspace{-2mm} \small\sl
Laboratoire Jean Kuntzmann / Universit\'e J. Fourier\\
\vspace{-2mm} \small\sl BP 53, 38041 Grenoble Cedex 9, France}

\date{}
\maketitle

\begin{abstract} Filtering problems  with general exponential quadratic criteria  are investigated for Gauss-Markov processes.
In this setting, the  Linear Exponential Gaussian  and Risk-Sensitive  filtering problems are solved
and it is shown that they
 may have different solutions.
\end{abstract}

\vspace*{4mm}


\vspace*{0.25cm}


\date{}

\maketitle

\section{Introduction}
The so-called linear exponential Gaussian (LEG) and risk-sensitive (RS) filtering problems involve criteria which are exponentials of integral cost  functionals.  Before our paper \cite{AMM}, numerous results had been already reported in specific
models, specially around  Markov models, but
without exhibiting  the relationship between these two problems.
See, \textit{e.g.}, Whittle \cite{whittle1}, Speyer  \textit{et
al.} \cite{speyer}, Elliott \textit{et al}. \cite{elliott4}, \cite{elliott3} and
\cite{elliott1}
for contributions.
 In our paper \cite {AMM}, we have  solved the LEG and RS filtering problems for general
 Gaussian processes in the particular setting where the  functional in the exponential is a \textit{singular} quadratic functional. Moreover we have  proved that actually in this case the solutions coincide.
 In the present  paper the problems are revisited for Gauss-Markov processes but with a \textit{nonsingular} quadratic functional in the exponential.
 In this setting  the solutions are exhibited and we propose an example to show that they may be different.


It what follows all random variables and processes are defined on
a given stochastic basis ($\Omega,{\cal F},({\cal
F}_t),\space\Pg$) satisfying the usual conditions and processes
are $({\cal F}_t)$-adapted.
We deal
with a signal process $X=(X_t,\,t\geq 0)$ in $\Rg$ governed by the linear equation
\begin{equation}\label{eq:model}
dX_t= a_{t} X_t \, dt + dB_t,\, X_0=0,
\end{equation}
and an observation process $Y=(Y_t,\,t\geq 0)$ in $\Rg$
governed by the linear equation
\begin{equation}
\label{eq:observations}
dY_t=A_{t}X_t\,dt+d\widetilde{B}_{t}\,,Y_0=0,\quad t\geq 0\,.
\end{equation}
Here $a=(a_{t}, t\geq 0)$ and $A=(A_{t}, t\geq 0)$ are continuous real-valued  deterministic functions,
$B=(B_t,\, t\geq 0)$ and $\widetilde{B}=(\widetilde{B}_t,\, t\geq 0)$
are  independent
standard  one dimensional Brownian motions. Clearly the pair $(X,Y)$ is Gaussian.\\
For  a given  continuous deterministic function  $\Lambda=(\Lambda_s,\,0\le s\le T)$ with values in  the set of nonnegative definite symmetric
$2\times 2$ matrices
$$
\Lambda_s =\left(\begin{smallmatrix}
\Lambda_{11}(s) & \Lambda_{12}(s)
\\
\Lambda_{12}(s) & \Lambda_{22}(s)
\end{smallmatrix}\right),
$$
such that  $\Lambda_{22}(s)\neq 0,$ let us define by $\bar{h}\in {\cal H}$  the solution of the
\textit{LEG type filtering problem }:

\begin{equation}\label{defrssex}
\bar{h}=\argmin_{h\in\mathcal{H}} \Eg\Big[\mu \exp \left\{\frac{\mu}{2}\int_0^T (X_s \, h_s)\Lambda_{s} \left(X_s \atop h_s \right) \, ds \right\} \Big].
\end{equation}
In this definition $\mu$ is a real parameter and  $h=(h_{s},0\le s\le T)\in {\cal H}$ means that $h$ is a
$({\cal Y}_s)$-adapted continuous process  where
$({\cal Y}_s)$ is the natural filtration of $Y$, \textit{i.e.},
${\cal Y}_s=\sigma(\{Y_u\, ,\, 0\leq u\leq s\}),0\le s\le T$.

We can also define
$\hat{h}$  as a solution of
the following recursive equation, which is the basic definition of the \textit{ RS type filtering problem}:
\begin{equation}\label{rdeex}
\widehat{h}_t= \argmin_{g\in {\cal Y}_t} \Eg\Big[\left.\mu\exp\left\{\frac{\mu}{2}(X_t \, g) \Lambda_{t} \left(X_t \atop g\right) +  \frac{\mu}{2} \int_0^t (X_s \, \widehat{h}_s) \Lambda_{s} \left(X_s \atop \widehat{h}_s \right) \, ds \right\} \right/ {\cal Y}_t\Big],
\end{equation}
where  $g\in {\cal Y}_t$ means that $g$ is a ${\cal Y}_t$-measurable variable.

It is clear that \textit{risk-neutral} versions of these two problems (namely, dropping the exponentials in definitions \eqref{defrssex}-\eqref{rdeex}, \textit{i.e.}, simply with quadratic criteria)  are ``equivalent'':
$$
\bar{h}_t=\widehat{h}_t=- \frac{\Lambda_{12}(t)}{\Lambda_{22}(t)}\cdot \pi_t(X),
$$
where for any process $\eta=(\eta_t\,,\, t\in [0,T])$ such that
$\Eg |\eta_t|<+\infty$, the notation
$\pi_t(\eta)$ is used for the conditional expectation of $\eta_t$ given the
$\sigma$-field
${\cal
Y}_t,$
$$
\pi_t(\eta) = \Eg(\eta_t/{\cal Y}_t)\,.
$$

One question that we want to discuss in this paper is the possible ``equivalence'' of the problems \eqref{defrssex} and \eqref{rdeex}.
In our paper \cite{AMM}, we have proved that when the  quadratic functional involved in the exponential is \textit{singular}, namely when  matrices
$\Lambda_{s}$ are singular,\textit{ i.e.},
$\Lambda_{11}=\Lambda_{22}=-\Lambda_{12}$, the equality
$\bar{h}=\hat{h}$ holds, even in a non Markovian setting. Here below
a simple example where
$
\bar{h}\neq \wh
$
is proposed which  shows that if  the  quadratic functional  is \textit{nonsingular}
then the answer may be negative even for the Markovian  model~\eqref{eq:model}--\eqref{eq:observations}.\\
The paper is organized as follows. Preparing for the analysis of the filtering problems, in Section \ref{CMNS}  a
 Cameron-Martin type formula for the {\em conditional
Laplace transform} of a quadratic functional of the involved signal
process is derived. Then in  Section \ref{solfpbs}  the LEG and RS filtering problems in the nonsingular setting
are solved.
Finally, Section \ref{EX} is devoted to the analysis of the announced example which shows the discrepancy between the two filtering problems.

\section{Conditional version of a Cameron-Martin formula}\label{CMNS}
Actually, the resolution of the LEG and RS filtering problems is based on
 a conditional version of a Cameron-Martin formula and we follow the same lines as in our paper \cite{AMM}.\\
In the present Section, the process  $X=(X_t, t\geq 0)$ is an arbitrary continuous Gaussian process with mean function
$m=(m_t, t\geq 0)$ and covariance function  $K =(K(t,s), t\geq 0,
s\geq 0)$, \textit{i.e.},
$$
\Eg X_t=m_t,\quad\Eg (X_t-m_t)(X_s-m_s)=K(t,s)\,, \quad
t\geq 0\,,\; s\geq 0\,.
$$
We are interested in the explicit representation of
\begin{multline}
{\cal I}_{T}= \Eg\Big[\left.\mu\exp\left\{\frac{ \mu}{2}(X_T \, g) M \left(X_T \atop g\right) +
\frac{ \mu}{2} \int_0^T (X_s \, h_s) \Lambda_{s} \left(X_s \atop h_s \right) \, ds \right\} \right/ {\cal Y}_T\Big],\\
\end{multline}
for any variable $g\in {\cal Y}_{_{T}}$ and process $ h\in\mathcal{H}$,
and with symmetric deterministic nonnegative definite
matrices $M$ and $\Lambda_{s}$.

Let us formulate the
 condition  $(C_{\mu})$:
\begin{enumerate}
  \item[$(C_{\mu})$] the Riccati-Volterra  equation

\begin{equation}\label{eq:gammagen}
\bar{\gamma}(t,s)=K(t,s) - \int_0^s
\bar{\gamma}(t,r) [A_{r}^{2}-\mu \Lambda_{11}(r)]
\bar{\gamma}(s,r) dr,\, 0\le s \le t \le T,
\end{equation}
 has a unique and bounded solution on $\{(t,s):0\le s \le t \le T\}$,  such that $\bar{\gamma}(t,t)\geq 0$ for $0\leq t\leq T$
 and moreover
$$ \displaystyle{1 -\mu
 M_{11}\bar{\gamma}(T,T)}> 0.
 $$
\end{enumerate}
Notice that for all  $\mu$ {\em negative} the condition $(C_{\mu})$ is satisfied
 and if $\mu $ is  {\em positive}, the condition $(C_{\mu})$ is satisfied for $\mu$  sufficiently small, for example, those such  that for any
  $ t \le T \quad  A_{t}^{2}-\mu \Lambda_{11}(t)$ is nonnegative (\textit{cf.} Lemma 2 \cite{AMM}).

Now we claim the following extension of the $1-D$ version of Proposition 2 \cite{AMM}:

\begin{prop}\label{C-M leg}
Suppose that the condition $(C_{\mu})$ is satisfied. Let  $Z^{h}=(Z^{h}_{s},
0\le s\le T)$  be the unique solution of the It\^o-Volterra equation
\begin{equation}\label{eqzhgen}
Z^{h}_{t} = m_{t}+ \displaystyle{
 \int_{0}^{t}\bar{\gamma}(t,s)\mu [\Lambda_{11}(s)Z^{h}_{s}+\Lambda_{12}(s)h_{s}]ds +
 \int_{0}^{t}\bar{\gamma}(t,s)A_{s}[dY_{s} - A_{s}Z^{h}_{s}ds]},
  \end{equation}
and $\bg_{_{XX}}(t)= \bar{\gamma}(t,t), 0\le t\le T$ where $\bar{\gamma}$ is the unique solution of equation \eqref{eq:gammagen}.
Then the following equality holds:
\begin{multline}\label{eq:claim}
 {\cal I}_{T}= (1-\mu M_{11}{\bg_{_{XX}}}(T))^{-1/2} \exp\left\{\frac{\mu}{2}\int_0^T {\bg_{_{XX}}}(s) \Lambda_{11}(s) \, ds \right\} \times \\
\times \exp \left\{ \frac{ \mu}{2} (Z_T^{h} \, g) G_T \left(Z_T^{h} \atop g\right) +\frac{ \mu}{2}\int_0^T (Z_s^{h} \, h_s) \Lambda_s \left(Z_s^{h} \atop h_s \right) \, ds \right\} \times
\\
\times \exp \left\{ \int_0^T A_s (Z_s^{h}-\pi_s(X)) \, d\nu_{s} - \frac{1}{2} \int_0^T |A_s(Z_s^{h} - \pi_s(X))|^2 \, ds \right\},
\end{multline}
where
\begin{equation}\label{def:G}
G_{T} = (1-\mu M_{11} {\bg_{_{XX}}}(T))^{-1} \left(\begin{array}{cc}
M_{11} & M_{12} \\
M_{12} & M_{22} -\mu{\bg_{_{XX}}}(T) \mathop{det}(M)
\end{array}\right),
\end{equation}
and $(\nu_{t}, \, t\ge 0)$ is the innovation process associated to  $Y$,\textit{ i.e.},

\begin{equation}\label{def nu}
    d\nu_{t}= dY_t - A_t\pi_t(X)dt,\, \nu_{0}=0.
\end{equation}

\end{prop}
\begin{rem}
\begin{enumerate}
  \item
  Note that in the singular case where $ M_{11}=M_{22}=-M_{12}$ and $\Lambda_{11}=\Lambda_{22}=-\Lambda_{12}$
  Proposition \ref{C-M leg} reduces to the $1-D$ version of Proposition 2 \cite{AMM}.

  \item
  Note also that the condition $(C_{\mu})$ implies that $G_T$ is nonnegative definite.

\end{enumerate}

\end{rem}
\paragraph{Proof of the Proposition \ref{C-M leg}} The proof is based on the ideas developed in the
proof of  Propositions 1 and 2 \cite{AMM}. Actually, it is sufficient to work with $\mu < 0$
since  the result will be  valid
for sufficiently small $\mu > 0$  because of the analytical properties of the
involved functions.
To simplify the notations we work with $\mu = -1$; then for the
general situation it is sufficient to replace $M$ and $\Lambda$ by $-\mu M$ and $-\mu \Lambda$ respectively.

Let us introduce the auxiliary observations $(\bar{Y}_{t},\, 0 \le t \le T)$ such that:
\begin{equation}\label{eq:Ybar}
\left\{
\begin{array}{lcl}
d\bar{Y}_t^1 & = & dY_t,
\\
d\bar{Y}_t^{2,3} & = & d\bar{B}_t + \Lambda_t^{\frac{1}{2}} \left(X_t \atop h_t \right) \, dt,
\end{array}
\right.
\end{equation}
where $\bar{B}=(\bar{B}_t, t\geq 0)$
denotes a 2-D
standard  Brownian motion, independent of ($X,\,\widetilde{B}$).\\
Below, for any process $\eta=(\eta_t\,,\, t\in [0,T])$ such that
$\Eg |\eta_t|<+\infty$,  the notation
$\bar{\pi}_t(\eta)$ is used for the conditional expectation of $\eta_t$ given the \textit{auxiliary}
$\sigma$-field
${\cal
\bar{Y}}_t=\sigma(\{\bar{Y}_s\, ,\, 0\leq s\leq t\}),\,
\bar{\pi}_t(\eta) = \Eg(\eta_t/{\cal\bar{ Y}}_t)\,.$\\
Let also $\xi_{t}$ be defined by
\begin{equation}\label{eq xi}
d\xi_t=(X_t \, h_t) \Lambda_t^{\frac{1}{2}} d\bar{Y}^{2,3}_{t},\, \xi_0=0.
\end{equation}

We see that the conditional  distribution of $(X_{t},\xi_{t})$ given
${\cal\overline{ Y}}_t$ is Gaussian with the
conditional expectation  $\left( \bar{\pi}_{t}(X),\,  \bar{\pi}_{t}(\xi) \right)$
and the conditional covariance
$\left(
                               \begin{array}{cc}
                                \bar{\gamma}_{_{XX}}(t)  & \bar{\gamma}_{X\xi}(t) \\
                                 \bar{\gamma}_{X\xi}(t) & \bar{\gamma}_{\xi\xi}(t)
                               \end{array}
                             \right),$
where
$$
\bar{\gamma}_{_{XX}}(t)=\Eg[(X_t-\bar{\pi}_t(X))^{2}/{\cal \overline{Y}}_t],
$$
\begin{equation}\label{def gamma Xxi}
\bar{\gamma}_{_{X\xi}}(t)=\Eg[(X_t-\bar{\pi}_t(X))(\xi_t-\bar{\pi}_t(\xi))/{\cal\overline{ Y}}_t]\,,
\end{equation}
and
\begin{equation}\label{def gamma xixi}
\bar{\gamma}_{_{\xi\xi}}(t)=\Eg[(\xi_t-\bar{\pi}_t(\xi))^{2}/{\cal\overline{ Y}}_t]\,.
\end{equation}
Proceeding as in \cite{AMM} Section 2.2, we obtain that
\begin{itemize}
  \item the conditional variance
$\bar{\gamma}_{_{XX}}(t)$ is deterministic
and actually  nothing but the variance of the filtering error,\textit{ i.e.},
\begin{equation}\label{def gamma XX}
\bar{\gamma}_{_{XX}}(t)=\Eg[(X_t-\bar{\pi}_t(X))^{2}]\,,
\end{equation}
given by $\bar{\gamma}_{_{XX}}(t)= \bar{\gamma}(t,t)$, where $\bar{\gamma}(t,s)$ is the unique solution of the equation \eqref{eq:gammagen} with $\mu=-1$,
  \item the difference
\begin{equation}\label{Zhdef}
Z^{h}_{t}= \bar{\pi}_{t}(X)-\bar{\gamma}_{_{X \xi}}(t)
\end{equation}
is ${\cal Y}_t$-measurable  and is the unique solution of the equation~\eqref{eqzhgen} with $\mu=-1$,
  \item  $\bar{\pi}_t(\xi)$ is the solution of the equation:
\begin {equation}\label{pi s xi tc}
\begin{array}{ccl}
\bar{\pi}_t(\xi)&=&\displaystyle{\int_{0}^{t}\bar{\gamma}_{_{XX}}(s)
\Lambda_{11}(s) ds + \int_0^t(\bar{\pi}_s(X),h_{s})\Lambda_{s}\left(\bar{\pi}_s(X) \atop h_{s}\right)ds}\\
&+&\displaystyle{\int_0^t[(\bar{\pi}_s(X),h_{s})+\bar{\gamma}_{_{X\xi}}(s)(1,0)]
\Lambda_{s}^{\frac{1}{2}}d\bar{\nu}_s^{2,3}
+\int_0^t\bar{\gamma}_{_{X\xi}}(s)A_{s}d\bar{\nu}_s^{1}}\,,
\end{array}
\end{equation}
\item the conditional variance $\bar{\gamma}_{_{\xi\xi}}(t)$ satisfies the equation:
\begin{equation}\label{eq gamma xi}
\begin{array}{ccl}
\bar{\gamma}_{_{\xi\xi}}(t)&=&
\displaystyle{\int_{0}^{t}\bar{\gamma}_{_{XX}}(s)
\Lambda_{11}(s) ds + 2\int_{0}^{t}\bar{\gamma}_{_{X\xi}}(s)(1,0)\Lambda_{s}^{\frac{1}{2}}d\bar{\nu}_{s}^{2,3}}
\\
&-&
\displaystyle{\int_{0}^{t}\bar{\gamma}_{_{X\xi}}(s)[\Lambda_{11}(s)
+
A^{2}_{s}]\bar{\gamma}_{_{X\xi}}(s) ds }
\\&+&
2\displaystyle{ \int_{0}^{t}(\bar{\pi}_s(X),h_{s})
\Lambda_{s}\left(1 \atop 0\right) \bar{\gamma}_{_{X\xi}}(s)ds}\,,
\end{array}
\end{equation}
\end{itemize}
where $\bar{\nu}_{t}=(\bar{\nu}_{t}^{1}, [\bar{\nu}_{t}^{2,3}]^{\prime})^{\prime}$ is the $3-D$ innovation process associated to the auxiliary observations $\bar{Y}$, \textit{ i.e.},

\begin{equation}\label{def barnu}
    d\bar{\nu}_{t}= d\bar{Y}_t - \left(
                               \begin{array}{ccc}
                                A_{t}  & & 0 \\
                                 &\Lambda_{t}^{\frac{1}{2}} &
                               \end{array}
                             \right)\left(\bar{\pi}_{t}(X) \atop h_{t}\right)dt,\, \bar{\nu}_{0}=0.
\end{equation}

Now we turn to the proof of equality \eqref{eq:claim} for $ \mu =
-1$.\\
Let $\rho_{t}$ be defined by:
\begin{equation}\label{eq ro}
\rho_{t}=\exp\left\{-\int_0^t
(X_{s},h_{s})\Lambda^{\frac{1}{2}}_{s}\,d\bar{B}_s -
\frac{1}{2}\int_0^t (X_{s},h_{s})\Lambda_{s}\left(X_{s} \atop h_{s}\right) ds
\right\}.
\end{equation}

At first we  note that the same arguments that we have used in the proof of the
Proposition 1 \cite{AMM}  give  the equality:

\begin{equation}\label{bayesfor}
{\cal I}_{T} =
\frac{\Eg[\exp(-\frac{1}{2}(X_T \, g) M \left(X_T \atop g\right)-\xi_{T})/
{\cal \overline{Y}}_T]}{\Eg[\rho_{T}/
{\cal \overline{Y}}_T]}
=\frac{\varphi_{1}(T)}{\bar{\pi}_{T}(\rho)}\,,
\end{equation}
where $\xi_{T}$ and $\rho_{_{T}}$  are defined by (\ref{eq xi}) and (\ref{eq ro})
 respectively. But the conditional Gaussian
properties of the pair $(X,\xi)$ given    ${\cal \bar{Y}}_T $ gives  the
following (see for example \cite {lipshi1}, Lemma 11.6):
\begin{equation}\label{numer}
\begin{array}{ccl}
 \ln  \varphi_{1}(T)&=& -\displaystyle{\frac{1}{2} \ln (1 +
 M_{11}\bar{\gamma}_{_{XX}}(T))} \\
 &-& \displaystyle{\frac{1}{2}(Z^{h}_{T},g)G_T\left(Z_T^{h} \atop g\right)}
-\displaystyle{ \bar{\pi}_{T}(\xi)
+\frac{1}{2}\bar{\gamma}_{_{\xi \xi}}(T)}\,,
\end{array}
\end{equation}
where the terms $G_{T},\,\bar{\gamma}_{_{\xi \xi}}(T),\,\bar{\gamma}_{_{XX}}(T)$ and $ Z^h_{T}$  are defined by the equations \eqref{def:G},
\eqref{def gamma xixi}, \eqref{def gamma XX} and \eqref{Zhdef}   respectively.

Now it follows from (\ref{numer}) that to prove the statement of the Proposition it is
sufficient to write the expression for
$$
\Psi_{_{T}}=\frac{\exp (- \bar{\pi}_{T}(\xi)
+\frac{1}{2}\bar{\gamma}_{_{\xi \xi}}(T))}{\bar{\pi}_{_{T}}(\rho)}.
$$

Since
$$
d\rho_{t}=-\rho_{t}(X_t \, h_t)\Lambda_{t}^{\frac{1}{2}}d\bar{B}_{t}\,,\rho_{0}=1\,,
$$
thanks to the general filtering theorem \cite[Theorem 7.16]{lipshi1} we can write
\begin{equation}\label{denumer}
\begin{array}{ccl}
d\bar{\pi}_{t}(\rho)&=&\displaystyle{\bar{\pi}_{t}(\rho)
\left\{-(\pi_{t}(X),h_{t})\Lambda_{t}^{\frac{1}{2}}d\bar{\nu}_{t}^{2,3}\right.}\\
&+&
\displaystyle{\left.\left[\frac{\bar{\pi}_{t}(\rho X)}{\bar{\pi}_{t}(\rho)}
-\bar{\pi}_{t}(X)\right]A_{t}d\bar{\nu}_{t}^{1}\right\}}\,.
\end{array}
\end{equation}

We note that the
classical Bayes formula gives that
$$
\frac{\bar{\pi}_{t}(\rho X)}{\bar{\pi}_{t}(\rho)}
=\pi_{t}(X).
$$
Hence
$$
d\bar{\pi}_{t}(\rho)=\displaystyle{\bar{\pi}_{t}(\rho)
\left\{-(\bar{\pi}_{t}(X),h)\Lambda_{t}^{\frac{1}{2}}d\bar{\nu}_{t}^{2,3}
+ \left[\pi_{t}(X)
-\bar{\pi}_{t}(X)\right]A_{t}d\bar{\nu}_{t}^{1}\right\}}\,,
$$
or, equivalently:
\begin{equation}\label{denumerc}
\begin{array}{ccl}
\bar{\pi}_{_{T}}(\rho)&=&\displaystyle{\exp \left\{-\int_{0}^{T}(\bar{\pi}_{s}(X),h_{s})\Lambda_{s}^{\frac{1}{2}}
d\bar{\nu}_{s}^{2,3}
+\int_{0}^{T}\left[\pi_{s}(X)
-\bar{\pi}_{s}(X)\right]A_{s}d\bar{\nu}_{s}^{1}\right.}\\
&-&\displaystyle{\frac{1}{2}\left.\int_{0}^{T}|A_s\left[\pi_{s}(X)
-\bar{\pi}_{s}(X)\right]|^{2}ds
-\frac{1}{2}\int_{0}^{T}\|(\bar{\pi}_{s}(X),h_{s})\Lambda_{s}^{\frac{1}{2}}\|^{2}ds \right\}.}\\
\end{array}
\end{equation}

The equalities (\ref{pi s xi tc}), (\ref{eq
gamma xi}) and (\ref{denumerc}) imply:
$$
\ln(\Psi_{_{t}})=-\frac{1}{2}\int_{0}^{T}
 \bar{\gamma}(s,s)\Lambda_{11}(s)ds -
 \frac{1}{2}\int_{0}^{T}(\bar{\pi}_s(X),h_{s})\Lambda_{s}\left(\bar{\pi}_{s}(X)\atop h_{s}\right)ds
$$
$$
+\int_{0}^{T}(\bar{\pi}_s(X), h_{s})
\Lambda_{s}\left(1 \atop 0\right)\bar{\gamma}_{_{X\xi}}(s)ds -\frac{1}{2} \int_{0}^{t}\bar{\gamma}_{_{X\xi}}(s)\Lambda_{11}(s)
\bar{\gamma}_{_{X\xi}}(s) ds
$$
$$
+\int_0^T \left[\pi_{s}(X)
-\bar{\pi}_{s}(X)-\bar{\gamma}_{_{X\xi}}(s)\right]A_{s}d\bar{\nu}_s^{1}
$$
$$
-\frac{1}{2} \int_{0}^{T}\bar{\gamma}_{_{X\xi}}(s)A^{2}
\bar{\gamma}_{_{X\xi}}(s) ds +\frac{1}{2}\int_{0}^{T}|A_s\left[\pi_{s}(X)
-\bar{\pi}_{s}(X)\right]|^{2}ds.
$$
Replacing $d\bar{\nu}_t^{1}$ by $d\bar{\nu}_t^{1}=d\nu_t + A_{t}\left[\pi_{t}(X)
-\bar{\pi}_{t}(X)\right]dt $
we obtain that :
$$
\ln(\Psi_{_{t}})=-\frac{1}{2}\int_{0}^{T}
\bar{\gamma}(s,s)\Lambda_{11}(s)ds -
$$
$$
 -\frac{1}{2}\int_{0}^{T}(\bar{\pi}_s(X)-\bar{\gamma}_{_{X\xi}}(s),h_{s})
 \Lambda_{s}\left(\bar{\pi}_s(X) - \bar{\gamma}_{_{X\xi}}(s) \atop h_{s}\right)ds
$$
$$
+\int_0^T \left[\pi_{s}(X)
-\bar{\pi}_{s}(X)-\bar{\gamma}_{_{X\xi}}(s)\right]A_{s}d\nu_s
-\frac{1}{2} \int_{0}^{T}|A_{s}(\pi_{s}(X)
-\bar{\pi}_{s}(X)-\bar{\gamma}_{_{X\xi}}(s))|^{2}ds\,,
$$
and it gives the statement of the Proposition.
\begin{rem}
\begin{enumerate}
  \item Let us observe that for $\mu$ {negative} the proof of  Proposition~1 clarifies the probabilistic interpretation  of the ingredients  $\bg_{_{XX}}(t)$ and $Z^{h}_{t} $
in terms of an auxiliary risk-neutral filtering problem. They are nothing else but the  filtering error $\bg_{_{XX}}(t)$ (see equation \eqref{def gamma XX})
and  the difference
$Z^{h}_{t}= \bar{\pi}_{t}(X)-
\bar{\gamma}_{_{X \xi}}(t)$ (see equation \eqref{Zhdef}). It is worth mentioning that  $\bar{\pi}_t(X)$ and
      $\bar{\gamma}_{_{X\xi}}(t)$ are only  ${\cal\bar{ Y}}_t$-measurable variables but that the difference
$z^{h}_{t}=\bar{\pi}_t(X)-\bar{\gamma}_{_{X\xi}}(t)$  is actually a  ${\cal Y}_t$-measurable variable.

  \item If $\mu $  is  {positive}, but sufficiently small  in order that  the condition $(C_{\mu})$ is satisfied, due to  analytical properties of
 involved functions with respect to $\mu$,  equality \eqref{eq:claim} is still valid.
But it is worth  emphasizing that there is no  connection anymore between functions $Z^{h}$ and $\bg_{_{XX}}$
and a  risk-neutral filtering problem.

\end{enumerate}

\end{rem}
\section{Solution of the  filtering problems with exponentials of integral functionals  criteria}\label{solfpbs}

 Actually, in the particular Markov model \eqref{eq:model}--\eqref{eq:observations}, the equations \eqref{eq:gammagen}--\eqref{eqzhgen} can be transformed. Indeed, due to the specific structure of the covariance function $K$ of the signal process $X$, the solution of equation \eqref{eq:gammagen} is obtained in the form  $\bar{\gamma}(t,s)=\Pi_t\Pi_s^{-1}\bar{\gamma}_{_{XX}}(s), 0 \leq s\leq t$, where
$\Pi_s$ is the solution of the
differential equation $\dot\Pi_s=a(s)\Pi_s\,,\; s\geq 0\,,\Pi_0=1$ and $\bar{\gamma}_{_{XX}}(s)$ satisfies  the following differential equation:
\begin{equation}\label{eq:gammabis}
\dot{\bg}_{_{XX}}= 2a  \bg_{_{XX}} + 1 - \bg_{_{XX}}^2 [A^2 -\mu\Lambda_{11}],\, \bg_{_{XX}}(0) =0.
\end{equation}
Moreover, this particular form of $\bar{\gamma}(t,s)$ leads also to a differential equation for the solution $Z^{h}$ of \eqref{eqzhgen}:

\begin{equation}\label{eq:zh}
dZ_t^h=[a-\bar{\gamma}_{_{XX}} ( A^2-\mu\Lambda_{11} )] Z_t^h \, dt +\mu \bar{\gamma}_{_{XX}} \Lambda_{12} h_t \, dt + \bar{\gamma}_{_{XX}} A \, dY_t,\, Z_0^h=0.
\end{equation}

\subsection{Solution  of the LEG filtering problem}\label{LEGNS}

Let us formulate the following
 condition  $(C_{\mu}^{*})$:
\begin{enumerate}
  \item[$(C_{\mu}^{*})$] the forward and backward Riccati equations:
\begin{equation}\label{eq:gamma}
\dot{\bg}_{_{XX}}= 2a  \bg_{_{XX}} + 1 - \bg_{_{XX}}^2 [A^2 -\mu\Lambda_{11}],\, \bg_{_{XX}}(0)=0,
\end{equation}
\begin{equation}\label{eq:Riccati}
\dot{\Gamma} = -\frac{\mathop{det}(\Lambda)}{\Lambda_{22}} -2(a + \mu  \bg_{_{XX}}  \frac{\mathop{det}(\Lambda)}{\Lambda_{22}}) \Gamma -\mu \Gamma^2 \bg_{_{XX}}^2 [A^2 -\mu \frac{\Lambda_{12}^{2}}{\Lambda_{22}} ] ,\,\Gamma(T,T)=0,
\end{equation}
have  unique nonnegative and  bounded solutions  $ (\bg_{_{XX}}(t), \, 0\le t \le T) $ and $(\Gamma(T,t),  \, 0\le t \le T) $.
\end{enumerate}

Notice that for all  $\mu$ {\em negative} the condition $(C_{\mu}^{*})$ is satisfied
 and if $\mu $ is  {\em positive}, it is satisfied for $\mu$  sufficiently small.
\begin{prop}\label{LEGGM}
Suppose that the condition $(C_{\mu}^{*})$ is
satisfied. Let $\bar{h}=(\bar{h}_{t},\,0\le t \le T)$ such that:
\begin{equation}\label{barh:reprf}
\bar{h}_t = -\frac{\Lambda_{12}(t)}{\Lambda_{22}(t)}(1+\mu\bg_{_{XX}}(t) \Gamma(T,t) ) \, Z_t^{\bar{h}},
\end{equation}
where
$\Gamma(T,\cdot)$  is the solution of the backward Riccati equation \eqref{eq:Riccati}
and $ Z^{\bar{h}}= (Z_{t}^{\bar{h}}, \, 0\le t \le T)$ is the solution of the following equation:

\begin{equation}\label{eq:zhbarf}
 dZ_t^{\bar{h}}=[a+\mu\frac{\bg_{_{XX}}}{\Lambda_{22}}(\mathop{det}(\Lambda) -\mu\Lambda_{12}^{2}\bg_{_{XX}}\Gamma)]  Z_t^{\bar{h}} \, dt +
 A\bg_{_{XX}} \,[ dY_t -A Z_t^{\bar{h}}\,dt].
\end{equation}

Then
 $ \bar{h}$ is the
solution of the LEG filtering problem \eqref{defrssex} and
moreover, the corresponding optimal risk is given by
$$
\Eg\Big[\mu \exp \left\{\frac{\mu}{2}\int_0^T (X_s \,\bar{ h}_s)\Lambda_{s} \left(X_s \atop\bar{ h}_{s} \right) \, ds \right\} \Big]
$$
$$
=\mu\exp\left\{ \frac{\mu}{2} \int_0^T \bar{\gamma}_{_{XX}}(s) \Lambda_{11}(s) \, ds  +\frac{\mu}{2} \int_0^T \Gamma(T,s)A_s^{2} \bar{\gamma}_{_{XX}}^{2}(s) \, ds  \right\}.
$$

\end{prop}
\paragraph{Proof}
Of course, since we assume that condition $(C_{\mu}^{*})$ is satisfied, condition  $(C_{\mu})$ with $M=0$ is also fulfilled. Then we can apply
 the Cameron-Martin formula~\eqref{eq:claim} with $M=0$ (and hence in particular $G_{T}=0$). It gives that
\begin{multline}\label{costrepres}
\Eg \mu \exp \left\{ \frac{\mu}{2}\int_0^T (X_s \, h_s)\Lambda_{s} \left(X_s \atop h_s \right) \, ds \right\} =
\exp\left\{ \frac{\mu}{2} \int_0^T \bar{\gamma}_{_{XX}}(s) \Lambda_{11}(s) \, ds  \right\} \times
\\
\times
\Eg\mu \exp\left\{ \frac{\mu}{2} \int_0^T (Z^h_s \, h_s) \Lambda_s \left(Z^h_s \atop h_s \right) \, ds \,+\right.
\\
+ \left. \int_0^T A_s (Z_s^{h}-\pi_s(X)) \, d\nu_{s} - \frac{1}{2} \int_0^T |A_s(Z_s^{h} - \pi_s(X))|^2 \, ds \right\},
\end{multline}
where  for arbitrary $h\in {\cal H}$ the process $Z^h$ is the solution of equation \eqref{eq:zh}.
To find the solution of LEG filtering problem we propose to follow  the ideas of   \cite{bensoussan} and \cite{whittle2}, developed for the LEG control problem. Let us apply the   It\^o  formula to $\mu\Gamma(T,t) (Z_t^{h})^{2} $, where $\Gamma(T,\cdot)$ and $Z^{h}$ are the solutions of the equations \eqref{eq:Riccati} and \eqref{eq:zh} respectively. We see that
\begin{multline}\label{costrepresbis}
\Eg \mu \exp \left\{ \frac{\mu}{2}\int_0^T (X_s \, h_s)\Lambda_{s} \left(X_s \atop h_s \right) \, ds \right\} =
\\
\exp\left\{ \frac{\mu}{2} \int_0^T \bar{\gamma}_{_{XX}}(s) \Lambda_{11}(s) \, ds
+ \frac{\mu}{2} \int_0^T \Gamma(T,s)A_{s}^{2}\bar{\gamma}_{_{XX}}^{2}(s) \, ds  \right\}
\\
\times\Eg \mu \exp \left\{  \int_0^T  A_s (Z_s^{h}-\pi_s(X)+\mu\Gamma(T,s)\bar{\gamma}_{_{XX}}(s) Z^h_s) \, d\nu_{s}\, -\right.
\\
\left.-\frac{1}{2} \int_0^T | A_s (Z_s^{h}-\pi_s(X)+\mu\Gamma(T,s)\bar{\gamma}_{_{XX}}(s) Z^h_s)|^2 \, ds \right\}
\\
\times \exp \left\{  \frac{\mu}{2}\int_0^T\Lambda_{22}(s)\Big[h_s +
\frac{\Lambda_{12}(s)}{\Lambda_{22}(s)}(1+\mu\bar{\gamma}_{_{XX}}(s)\Gamma(T,s)) \, Z_s^{h}\Big]^{2}\, ds \right\}.
\end{multline}
Proceeding as in the proof of Theorem 1 \cite{AMM} we see that
Equation \eqref{costrepresbis} implies that the cost function \eqref{defrssex} (see also \eqref{costrepres})
has a uniform lower bound which is attained for $\bar{h}$ defined by the equation \eqref{barh:reprf}.

\begin{rem}
\begin{enumerate}
  \item It is clear that in the {singular} case where
$
\Lambda_{11}=\Lambda_{22} = -\Lambda_{12},
$
equation \eqref{eq:Riccati}  implies that $\Gamma\equiv 0$ and therefore $Z^{\bar{h}}=\bar{h}$
(\textit{cf}.\cite{AMM}).
  \item But in the general case $\Gamma$ may depend  on a~$T$ and as a consequence, $\bar{h}_t$ may also depend on ~$T$.
An example of such a dependence  will be given below.
Of course, by its  definition $\wh_t$ does not depend on~$T$ and hence $\bar{h} \ne \wh$ in this example.

\end{enumerate}

\end{rem}

\subsection{Solution of the RS filtering problem.}\label{RSNS}
Let us formulate the following
 condition  $(C_{\mu}^{**})$:
\begin{enumerate}
  \item[$(C_{\mu}^{**})$]  the Riccati equation \eqref{eq:gamma} has a unique, nonnegative and  bounded solution on $[0,T]$
  such that for $0 \le t \le T$
$$ \displaystyle{1 -\mu \bar{\gamma}_{_{XX}}(t)
 \Lambda_{11}(t)}> 0.
 $$
\end{enumerate}

\begin{prop}\label{RSGM}
Suppose that the condition $(C_{\mu}^{**})$ is
satisfied. Let $\wh=(\wh_{t},\,0\le t \le T)$ such that:
\begin{equation}\label{wh:represent}
\wh_t = - \frac{\Lambda_{12}(t)}{\Lambda_{22}(t)} [1-\mu{\bg_{_{XX}}(t)} \Lambda_{22}^{-1}(t)  \mathop{det}(\Lambda_{t})]^{-1}\, Z_t^{\wh},
\end{equation}

where
$ Z_t^{\wh}=(Z_{t}^{\wh}, \, 0\le t \le T)$ is the solution of the following equation:

\begin{equation}\label{eq:zhbarf}
 dZ_t^{\wh}=[a+\mu \bg_{_{XX}} \mathop{det}(\Lambda)  \frac{1-\mu\bg_{_{XX}}\Lambda_{11}\mathop{det}(\Lambda)}{\Lambda_{22} -\mu\bg_{_{XX}}\mathop{det}(\Lambda)}]  Z_t^{\wh} \, dt +
 A\bg_{_{XX}} \,[ dY_t -A Z_t^{\wh}\,dt].
\end{equation}

Then $ \wh$ is the
solution of the RS filtering problem \eqref{rdeex}.

\end{prop}
\paragraph{Proof}
Again, since we assume that condition $(C_{\mu}^{**})$ is satisfied, for any fixed $t \le T$, we can apply
 the Cameron-Martin formula~\eqref{eq:claim} with $t$ in place of $T$ and $\Lambda_{t}$ in place of $M$. It gives that
\begin{multline}\label{eq:claimrs}
 {\cal I}_{t}= (1-\mu \Lambda_{11}(t){\bg_{_{XX}}}(t))^{-\frac{1}{2}} \exp\left\{\frac{\mu}{2}\int_0^t {\bg_{_{XX}}}(s) \Lambda_{11}(s) \, ds \right\} \times \\
\times \exp \left\{ \frac{ \mu}{2} (Z_t^{\wh} \, g) G_t \left(Z_t^{\wh} \atop g\right) +\frac{ \mu}{2}\int_0^t (Z_s^{\wh} \, \wh_s) \Lambda_s \left(Z_s^{\wh} \atop \wh_s \right) \, ds \right\} \times
\\
\times \exp \left\{ \int_0^t A_s (Z_s^{\wh}-\pi_s(X)) \, d\nu_{s} - \frac{1}{2} \int_0^t |A_s(Z_s^{\wh} - \pi_s(X))|^2 \, ds \right\},
\end{multline}
where
\begin{equation}
G_{t} = (1-\mu \Lambda_{11}(t) {\bg_{_{XX}}}(t))^{-1} \left(\begin{array}{cc}
\Lambda_{11}(t) & \Lambda_{12}(t) \\
\Lambda_{12}(t) & \Lambda_{22}(t) -\mu{\bg_{_{XX}}}(t) \mathop{det}(\Lambda_{t})
\end{array}\right).
\end{equation}
Since  $\wh_{s}, 0 \le s <t$  is fixed, the optimization of the quadratic form $ (Z_t^{\wh} \, g) G_t \left(Z_t^{\wh} \atop g\right)$ in the equality \eqref{eq:claimrs}  with respect to $g$   gives the value of $\wh_t$:
$
\wh_t= -\frac{G_t^{12}}{G_t^{22}} \, Z_t^{\wh},
$
or equivalently
\begin{equation}
\wh_t = - \frac{\Lambda_{12}(t)}{\Lambda_{22}(t)} [1-\mu{\bg_{_{XX}}(t)} \Lambda_{22}^{-1}(t)  \mathop{det}(\Lambda)]^{-1}\, Z_t^{\wh},
\end{equation}
where $ Z_t^{\wh}$ is defined by  the equation \eqref{eq:zh} with $h=\wh$ and so by the equation \eqref{eq:zhbarf}.

\begin{rem}
\begin{enumerate}
  \item It is clear that for the  {singular} case where
$
\Lambda_{11}=\Lambda_{22} = -\Lambda_{12},
$
equalities \eqref{wh:represent} and \eqref{eq:zhbarf} imply that  $\wh=Z^{\wh}=Z^{\bar{h}}=\bar{h}$ (\textit{cf}.\cite{AMM}).

  \item Let us emphasize  that, of course, $\wh_t$ does not depend on $T$ and so generally $\wh_t\neq \bar{h}_t$.

\end{enumerate}
\end{rem}
Of course the RS filtering problem can be solved for  an arbitrary continuous Gaussian process $X=(X_t, t\geq 0)$  with mean function
$m=(m_t, t\geq 0)$ and covariance function  $K =(K(t,s), t\geq 0,
s\geq 0)$. To complete this section we propose the following generalization of Proposition \ref{RSGM} which can be proved by the same way.

\begin{prop}\label{RSGG}
Suppose that equation \eqref{eq:gammagen} has a unique and bounded solution $\bar{\gamma}= (\bar{\gamma}(t,s), \, 0 \le s \le t \le  T )$ such that
$\bar{\gamma}(t,t) \ge 0 $  and $ 1- \mu \bar{\gamma}(t,t)\Lambda_{11}(t) > 0$ for $0\le t \le T$.\\
 Let $\wh=(\wh_{t},\,0\le t \le T)$ such that:
\begin{equation}
\wh_t = - \frac{\Lambda_{12}(t)}{\Lambda_{22}(t)} [1-\mu{\bg_{_{XX}}(t)} \Lambda_{22}^{-1}(t)  \mathop{det}(\Lambda_{t})]^{-1}\, Z_t^{\wh},
\end{equation}
where $\bg_{_{XX}}(t)= \bar{\gamma}(t,t), 0\le t\le T$ and
$ Z^{\wh}=(Z^{\wh}_{t},
0\le t\le T)$ is the unique solution of the It\^o-Volterra equation:
\begin{multline}
Z^{\wh}_{t} = m_{t}+ \displaystyle{
 \mu\int_{0}^{t}\bar{\gamma}(t,s)\mathop{det}(\Lambda_{s})[\frac{1-\Lambda_{11}(s)\bar{\gamma}(s,s)\mathop{det}(\Lambda_{s})}
 {\Lambda_{22}(s)-\mu\bar{\gamma}(s,s)\mathop{det}(\Lambda_{s})}Z^{\wh}_{s}]ds }
 \\+
 \displaystyle{\int_{0}^{t}\bar{\gamma}(t,s)A_{s}[dY_{s} - A_{s}Z^{\wh}_{s}ds]}.
\end{multline}
\\
Then $\wh$  is the solution of the RS filtering problem \eqref{rdeex}.

\end{prop}

\section{Discrepancy between LEG and RS filtering problems: an example }\label{EX}
To show the possible dependence of the solution of  the LEG filtering problem on $T$ and so the discrepancy between LEG and RS filtering problems
we propose to take
$$
\Lambda=\left( \begin{array}{cc}
2 & -1 \\
-1 & 1
\end{array}\right), \quad a=0, \quad A=1,\quad  \mu=-1.
$$
In this case  equations \eqref{eq:gamma}  and \eqref{eq:Riccati} reduce to the following:
$$
\dot{\bg}_{_{XX}}=1-3\bg^2_{_{XX}}, \quad \bg_{_{XX}}(0)=0,
$$
\begin{equation}\label{Riccex}
\dot{\Gamma}=2\bg_{_{XX}} \Gamma -1 + 2\bg^2_{_{XX}} \Gamma^2, \quad \Gamma(T,T)=0.
\end{equation}
Thus
$$
{\bg_{_{XX}}}(t)=\frac{1}{\sqrt{3}} \tanh \sqrt{3} t.
$$
Equation \eqref{Riccex} can be also solved explicitly using the classical linearization method for Riccati   equations:
$$
\Gamma=\varphi_1^{-1} \varphi_2,\quad \varphi_2(T)=0, \quad \varphi_1(T)=1,
$$
where
\begin{equation}\label{eqvarfi}
\left\{\begin{array}{lcl}
\dot{\varphi}_1 & = & -\bg _{_{XX}}\varphi_1 - 2\bg^2_{_{XX}} \varphi_2 ,\hfill  \\
\dot{\varphi}_2 & = & -\varphi_1 + {\bg_{_{XX}}} \varphi_2.
\end{array}\right.
\end{equation}
Hence
$$
\ddot{\varphi_2} = \varphi_2 (\dot{\bg}_{_{XX}} + 3 \bg^2_{_{XX}}) =\varphi_2,
$$
and thus
$$
\varphi_2(t) =  \sinh (T-t),
$$
$$
\varphi_1(t) = \cosh (T-t) + \frac{1}{\sqrt{3}} \tanh \sqrt{3}t \cdot \sinh(T-t),
$$
$$
\Gamma(T,t) = \frac{\sinh(T-t)}{\cosh(T-t) + \frac{1}{\sqrt{3}} \tanh \sqrt{3}t \cdot \sinh (T-t)}\,\cdot
$$
Representation \eqref{barh:reprf} gives that
$$
\bar{h}_t = (1-{\bg_{_{XX}}(t)}\Gamma(T,t))Z_t^{\bar{h}},
$$
or
$$
\bar{h}_t= \int_0^t \bar{H} (T,t,s) \, dY_s,
$$
where
$$
\bar{H}(T,t,s)= (1-{\bg_{_{XX}}}(t)\Gamma(T,t,s)) \Psi(T,t,s),
$$
with
$$
\left\{\begin{array}{l}
\Psi(T,s,s)= {\bg_{_{XX}}}(s), \\
\dot{\Psi} = (-2{\bg_{_{XX}}} -{\bg_{_{XX}}}^2 \Gamma) \Psi.
\end{array}\right.
$$
Equation~(\ref{eqvarfi}) gives that
$$
(\ln \Psi)^{\cdot} = \frac{1}{2} (\ln \varphi_1)^{\cdot} - \frac{3}{2} {\bg_{_{XX}}} = \frac{1}{2} (\ln \varphi_1)^{\cdot} - \frac{1}{2} (\ln \cosh \sqrt{3}t)^{\cdot},
$$
Finally
$$
\Psi(T,t,s)=c(s) \sqrt{\frac{\varphi_1(t)}{\cosh\sqrt{3}t}},
$$
and
$$
\bar{H}(T,t,s) = \frac{\sinh \sqrt{3}s \cdot \cosh (T-t)}{\sqrt{\alpha_t \alpha_s}},
$$
with
$$
\alpha_t = \frac{\sqrt{3}+1}{2} \cosh (T+(\sqrt{3}-1)t) + \frac{\sqrt{3}-1}{2} \cosh (T-(\sqrt{3}+1)t).
$$
The same calculations based on the equality \eqref{wh:represent} give the representation
$$
\wh_t= \int_0^t \widehat{H}(t,s) \, dY_s,
$$
where
$$
\widehat{H}(t,s) = \frac{1}{\sqrt{3}} \cdot \frac{(\cosh \sqrt{3}t)^{1/3} \sinh \sqrt{3}s}{(\cosh \sqrt{3}s)^{2/3}}\,\cdot
$$
\bibliographystyle{plain}

\end{document}